\documentclass[review,fleqn,12pt]{elsarticle}

\usepackage{tipa}
\usepackage{bbm}
\usepackage{amsfonts}
\usepackage{amssymb}
\usepackage[fleqn]{amsmath}
\usepackage{amsthm}
\usepackage{amsbsy}
\usepackage{appendix}
\usepackage{graphicx}
\usepackage{color}
\usepackage{amsmath} 
\allowdisplaybreaks[4]

\newtheorem{theorem}{Theorem}
\newtheorem{lemma}{Lemma}
\newtheorem{corollary}{Corollary}

\newdefinition{remark}{Remark}
\newdefinition{assumption}{Assumption}
\newdefinition{definition}{Definition}
\begin{document}
\begin{frontmatter}
\title{On Generalized Sub-Gaussian Canonical Processes and Their Applications}


\author {Yiming Chen}
\ead{chenyiming960212@mail.sdu.edu.cn}
\author{Yuxuan Wang \corref{cor1}} 
\ead{wangyuxuan@zju.edu.cn} %
\author{Kefan Zhu}
\ead{zhu_kefan@whu.edu.cn}
%
\address{Institute for Financial Studies,
                Shandong University,
                Jinan
                250100,
                PR China}
\address{School of Mathematical Sciences,
                Zhejiang University,
                Hangzhou
                310058,
                PR China}
\address{School of Mathematics and Statistics, Wuhan University, Wuhan 430072, PR China}
\cortext[cor1]{Corresponding author}

\begin{abstract}
We obtain the upper bound of tail probability of generalized sub-Gaussian canonical processes. It can be viewed as a variant of the Bernstein-type inequality in the i.i.d case, and we further get a tighter bound of concentration inequality through uniformly randomized techniques. A concentration inequality for general functions involving independent random variables is also derived as an extension. As for applications, we derive convergence results for principal component analysis and the Rademacher complexities method.
\end{abstract}

\end{frontmatter}
\section{Introduction}
Let $X=\left(X_1, X_2, \ldots\right)$ be an independent (centered) random variable sequence. Write
$$
Y_t=\sum_i t_i X_i \quad \text { for } t=\left(t_1, t_2, \ldots\right) \in \ell^2
$$

Our aim is to derive the upper bound of the tail probability of the canonical process $Y_{t}$ where $X_i$ is a more generalized random variable compared to the Gaussian case. Such a tail bound is always regarded as a useful tool linked with estimating the expected value of the supremum of $Y_{t}$, i.e., $\mathrm{E}\sup_{t \in T}Y_{t}$, which plays a crucial role in characterizing the structure of the space $T$, the most important case, centered Gaussian processes $\mathrm{E}\sup_{t \in T}\left(G_t\right)_{t \in T}$  has been well understood after being connected with the geometry of the metric space $(T,d)$. The study of $\mathrm{E}\sup_{t \in T}\left(G_t\right)_{t \in T}$ can be traced back to the works of the classical chaining methods, which was initially studied by Kolmogorov, and developed by, among others,  Dudley \cite{ddl}, Fernique \cite{XF}, Talagrand \cite{MT}, van Handel \cite{RVV}. Specifically, the well-known Fernique-Talagrand majorizing measures theorem states that,

$$
\frac{1}{L} \gamma_2(T, d) \leq \mathrm{E}\sup_{t \in T}\left(G_t\right)_{t \in T} \leq L \gamma_2(T, d),
$$
where
$$
\gamma_2(T, d):=\inf \sup _{t \in T} \sum_{n=0}^{\infty} 2^{n / 2} \Delta\left(A_n(t)\right).
$$

Here, $\left(\mathcal{A}_n\right)_{n \geq 0}$ is an increasing sequence that satisfies $\left|\mathcal{A}_n\right| \leq 2^{2^n}$. The infimum is taken over all admissible sequences, and $A_n(t)$ represents the unique set in $\mathcal{A}_n$ that contains $t$. Moreover, $\Delta(A)$ denotes the diameter of the set $A$ with the distance $d(t, s)=\left(\mathbb{E}\left(G_t-G_s\right)^2\right)^{1 / 2}$.

Considering that any centred separable Gaussian process has the form of $\left(\sum_{i=1}^{\infty} t_i g_i\right)_{t \in T}$ by the Karhunen-Loève representation \cite{MR}, where $ g_1, g_2,$ are i.i.d. Gaussian random variables. Consequently, there has been interest in exploring cases where $X_i$ is a more complex process. Bednorz and Latala \cite{BLA}, \cite{LTA} have made a significant contribution to solving the Bernoulli case. Talagrand \cite{MT}     consider the $X_i$ are independent random variables sequences with density proportional to $\mathrm{P}(X_i \geq t) \leq \exp (-|t|^\alpha)$,  where $\alpha \geq 1$, which is also an extension of the Gaussian case.


For the applications of the upper bound of the canonical processes, Mendelson \cite{M2010}, \cite{M2016} showed that it provides theoretical guarantees for the suprema of unbounded empirical processes. In the early investigation of Krahmer et al., \cite{KMR}, $\mathrm{E}\sup_{t \in T}Y_{t}$ is linked with the suprema of the chaos processes, which is closely related to the restricted isometry property of random matrices in compressed sensing. 

In this paper, we consider $X_i$ following a $\varphi$-sub-Gaussian distribution, a concept introduced by Kozachenko \cite{BK}. $\varphi$-sub-Gaussian distributions are widely applied in various fields, including the queueing theory\cite{KO1} and information theory \cite{KO}. We refer the reader to \cite{BK},\cite{CLLW} for a detailed exposition.  We will provide a thorough overview of the preliminary in the ensuing sections.


As an extension, we are concerned about the concentration of $Y_t$ when $X_i$ does not necessarily have zero mean in the i.i.d case, since limited studies focus on this topic where $X_i$ follows a $\varphi$-sub-Gaussian distribution. The key technique is the so-called uniformly randomized Markov's inequality introduced by Ramdas et al.\cite{ARM}. Briefly speaking, given a nonnegative random variable $X$ and $U \sim \operatorname{Unif}(0,1)$, then for any $a >0$, we have that 
$$
\mathrm{P}(X \geqslant U / a)=\mathrm{E}[\mathrm{P}(U \leqslant a X \mid X)]=\mathrm{E}[\min (a X, 1)].
$$

With the Uniformly-randomized skills, Ramdas et al.\cite{ARM} further obtained the Hoeffding inequality for sub-Gaussian random variables,
$$
\mathrm{P}\left(\bar{X}_n-\mathrm{E}[X] \geqslant \sigma \sqrt{\frac{2 \log (1 / \alpha)}{n}}+\sigma \frac{\log (U)}{\sqrt{2 n \log (1 / \alpha)}}\right) \leqslant \alpha ,
$$
where $\bar{X}_n$ is the mean of $\sigma$-sub-Gaussian random variables $X_1, \ldots, X_n$. Compared to the classical result,
$$
\mathrm{P}\left(\bar{X}_n-\mathrm{E}[X] \geqslant \sigma \sqrt{\frac{2 \log (1 / \alpha)}{n}}\right) \leqslant \alpha ,
$$
it is a sharper bound due to $\log U<0$.

Moreover, we consider the concentration inequality with a generalized function $f: \Omega^n \rightarrow \mathbb{R}$, rather than the form of a linear combination of independent random variables. Maurer \cite{M2012}, \cite{MP21} have given a similar result for both the sub-Gaussian and sub-exponential cases with the entropy method. Additionally, as stated in Klochkov and Zhivotovskiy \cite{KZ}, assumptions on function $f$ lead to a slightly stronger bound. We may refer to \cite{BL} for more details on entropy methods. 

  In our applications, we explore the extension of principal component analysis (PCA) in learning theory, and we also conduct research on Rademacher complexity methods, which are discussed in \cite{KP}, and for the classical sub-Gaussian case, further references can be found in \cite{BM} and \cite{MZ}.

The structure of the paper is as follows. The preliminary on $\varphi$-sub-Gaussian distributions are given in Section 2. Section 3 presents the main results and the proofs. Section 4 contains examples and applications that demonstrate our outcomes.

\subsection*{Notation}
Throughout this paper, $X=\left(X_1, \ldots, X_n\right)$ is a vector of independent random variable, $X^{*}=\left(X_1^{*}, \ldots, X_n^{*}\right)$  represents the independent copy of $X$. $c$ is denoted as a universal constant. If $H$ is a Hilbert space, then the Hilbert space of Hilbert-Schmidt operators consists of bounded operators $T$ on $H$ that satisfy $||T||_{HS}=\sqrt{\sum{ij}\left\langle Te_i, e_j\right\rangle_H^2}<\infty$. The inner product $\langle T, S\rangle_{HS}$ is defined as $\sum_{ij}\left\langle Te_i, e_j\right\rangle_H\left\langle Se_i, e_j\right\rangle_H$, where $\left(e_i\right)$ represents an orthonormal basis. For any $x \in H$, the operator $R_x$ is defined by $R_x y=\langle y, x\rangle x$, and it can be verified that $||R_x||_{HS}=||x||_H^2$.

\section{Preliminary}
In this section, we present definitions and properties of the $\varphi$-sub-Gaussian distribution, which serves as a valuable complement to sub-Gaussian processes in various applications. For instance, in situations where the summands exhibit large variances, potentially leading to a distribution that deviates from Gaussian behavior, the $\varphi$-sub-Gaussian property remains a reasonable assumption. Notably, all centered bounded processes possess $\varphi$-sub-Gaussian properties, and the sums of independent Gaussian processes and centered bounded processes also exhibit $\varphi$-sub-Gaussian characteristics.

\begin{definition}
A continuous, even, convex function $\varphi(x)$ defined for $x \in \mathbb{R}$ is referred to as an Orlicz N-function if it is monotonically increasing for $x>0, \varphi(0)=0, \varphi(x) / x \rightarrow 0$, when $x \rightarrow 0$, and $\varphi(x) / x \rightarrow \infty$, when $x \rightarrow \infty$.
\end{definition}

\begin{lemma}\label{book}(Buldygin and Kozachenko \cite{BK})
For any Orlicz $N$-function $\varphi$, the following results hold:
\begin{itemize}
\item for $\beta>1$, $\varphi(\beta x)\geq\beta\varphi(x)$ for $x\in\mathbb{R}$;
\item there exists a constant $c=c(\varphi)>0$, such that $\varphi(x)>cx$ for $x>1$;
\item the function $g(x)=\frac{\varphi(x)}{x}$ is monotone nondecreasing in $x$ for all $x>0$;
\item $\varphi(x)+\varphi(y)\leq\varphi(|x|+|y|)$ for $x,y\in\mathbb{R}$.
\end{itemize}
\end{lemma}

\begin{definition}
For any function $\varphi(x), x \in \mathbb{R}$, the function $\varphi^*(x):=\sup _{y \in \mathbb{R}}(x y-\varphi(y))$, $x \in \mathbb{R}$, is called the Young-Fenchel transform of $\varphi(\cdot)$, which is also known as convex conjugate.
\end{definition}


\begin{definition}
Given an Orlicz N-function $\varphi(\cdot)$, a zero-mean random variable $\xi$ is called $\varphi$-sub-Gaussian distribution if there exists a non-negative constant $a \geq 0$ such that the following inequality holds for all $\lambda \in \mathbb{R}$:
$$
\mathrm{E} \exp (\lambda \xi) \leq \exp (\varphi(\lambda a))
$$
\end{definition}


\begin{definition}\label{psi}
Let $\operatorname{Sub}_{\varphi}(\Omega)$ denote the space of $\varphi$-sub-Gaussian random variables, and then for any random variable $\xi \in \operatorname{Sub}_{\varphi}(\Omega)$,
$$
\mathrm{E} \exp \{\lambda \xi\} \leq \exp \left\{\varphi\left(\lambda \tau_{\varphi}(\xi)\right)\right\}, \quad \lambda \in \mathbb{R}
$$
holds, where
$$
\tau_{\varphi}(\xi):=\sup _{\lambda \neq 0} \frac{\varphi^{(-1)}(\log \mathrm{E} \exp \{\lambda \xi\})}{|\lambda|},
$$
and $\varphi^{(-1)}(\cdot)$ denotes the inverse function of $\varphi(\cdot)$.
\end{definition}

In fact, the norm $\tau_{\varphi}(\xi)$ is equivalent to 
$
\tau_{\varphi}(\xi)=\inf \{a \geq 0: E \exp \{\lambda \xi\} \leq \exp \{\varphi(a \lambda)\}, \lambda \in \mathrm{R}\}
$. It is evident that when $\varphi(\xi)=\xi^2 / 2$, the space is a  of sub-Gaussian variables space.

For the random variable $f(X)$, when considering the concentration of its expectation, the following definition regarding the fluctuations in the performance of the function $f(\cdot)$ for each component of $X$  is essential.

\begin{definition}
For a function $f: \Omega^n \rightarrow \mathbb{R}$, where $x=\left(x_1, \ldots, x_n\right) \in \Omega^n$ and $X=\left(X_1, \ldots, X_n\right)$ is a random vector with independent components in $\Omega^n$, the $k$-th centered conditional version of $f$ is the random variable that can be expressed as:
$$
f_k(X)(x)=f\left(x_1, \ldots, x_{k-1}, X_k, x_{k+1}, \ldots, x_n\right)-\mathrm{E}\left[f\left(x_1, \ldots, x_{k-1}, X_k^{\prime}, x_{k+1}, \ldots, x_n\right)\right] .
$$
\end{definition}

Observe that $f_k(X)$ is also a random-variable about $f_k(X): x \in \Omega^n \mapsto f_k(X)(x)$, which does not depend on the $k$-th coordinate of $x$.


\section{Main Results}\label{sec-Intro}
Our first result presents the tail bound of the $\varphi^*$-sub-Gaussian canonical process where each $X_i$ corresponds to a different distribution. To be more specific, for all $\lambda \in \mathbb{R}$, there exists a non-negative constant $a_i \geq 0$ such that $\mathrm{E} \exp(\lambda X_i)) \leq \exp(\varphi^*_i(\lambda a_i))$.  We state that the following result provides an upper bound under a non-identically distributed case for the $\varphi^*$-sub-Gaussian distribution. First, given $v>0$, for $t=\left(t_1, t_2, \ldots\right) \in \ell^2$, define
$$
\mathcal{N}_v(t)=\sup \left\{\sum_{i \geq 1} t_i b_i ; \sum_{i \geq 1} \varphi_i\left(b_i\right) \leq v\right\} .
$$
\begin{theorem}\label{2.3}
For $v>0, s \geq 1$, assume that $X_i$ is a sequence of independent $\varphi^*_i$-sub-Gaussian random variables with $\tau_{\varphi^*_i}(X_i)$ existing, the random variable $Y_{t}=\sum_{i \geq 1} t_i X_i$ satisfies 
$$\mathrm{P}\left(Y_{t} \geq 2 s K \mathcal{N}_v(t)\right) \leq \exp (-v s),$$

Specifically, let $X_i$ be independent and identically distributed $\varphi^*$-sub-Gaussian random variables,
$$
\mathrm{P}\left(Y_{t} \geq z\right) \leq \exp \left(-c\min\left( \varphi\left(\frac{z}{K_1 ||t||_1}\right),\frac{z^2}{K_2^2 ||t||_2^2}\right)\right),
$$
where $K=\displaystyle\max_{ i} \tau_{\varphi^*_i}(X_i)$, $K_1=\displaystyle\max_{ i} \tau_{\varphi^*}(X_i)$, $K_2=\max _i||X_i||_{\psi_1}$.
\end{theorem}



\begin{proof}[Proof of Theorem \ref{2.3}]
By the the convexity of the function $\varphi(x)$, for $s \geq 1$, we have that $\varphi_i\left(X_i / s\right) \leq \varphi_i\left(X_i\right) / s$, consider constants $b_i$ such that
$$
\sum_{i \geq 1} \varphi_i\left(b_i / s\right) \leq v.
$$
Recall that $\mathcal{N}_v(t)=\sup \left\{\sum_{i \geq 1} t_i b_i ; \sum_{i \geq 1} \varphi_i\left(b_i\right) \leq v\right\}$, so $
\sum_{i \geq 1} t_i b_i  \leq s \mathcal{N}_v(t)
$, then $s \mathcal{N}_v(t) \geq \mathcal{N}_{v s}(t)$. Assume $s=1$ without loss of generality, for $\lambda \geq 0$, based on the independence of random variables, we derive
$$
\begin{aligned}
\mathrm{E} \exp \left(\lambda Y_{t}\right)&=\mathrm{E} \exp \left(\lambda \sum_{i \geq 1} t_i X_i\right)=\prod_{i \geq 1} \mathrm{E} \exp \left(\lambda t_i X_i\right)\\
& \leq \prod_{i \geq 1} \exp \left(\varphi^*_i\left(\lambda t_i \tau_{\varphi^*_i}\left( X_i\right)\right)\right)\\
&= \exp \left(\sum_{i \geq 1} \varphi^*_i\left(\lambda t_i \tau_{\varphi^*_i}\left( X_i\right)\right)\right). \\
\end{aligned}
$$
Then, by Markov’s inequality, we have
$$
\begin{aligned}
\mathrm{P}\left\{Y_{t} \geq z\right\} & =\mathrm{P}\left\{\exp \left(\lambda Y_{t}\right) \geq \exp (\lambda z)\right\}  \leq e^{-\lambda z} \mathrm{E} \exp \left(\lambda Y_{t}\right)\\
&\leq\exp \left(-\lambda z +\sum_{i \geq 1} \varphi^*\left(\left(K \lambda\left|t_i\right|\right)\right)\right).\\
\end{aligned}
$$
Set $z=2 K \mathcal{N}_v(t)$, we control the second term inside the exponent in the above inequality, given numbers $x_i \geq 0$, for any $v\geq0$, if $\sum_{i \geq 1} \varphi_i\left(x_i\right) \leq u$, by the definition of $\mathcal{N}_v(t)$, we have $$\sum_{i \geq 1}\left|t_i\right| x_i \leq\mathcal{N}_v(t),$$ and since $\sum_{i \geq 1} \varphi_i\left(x_i\right) \geq 0$, we can conclude that
$$
\sum_{i \geq 1} \frac{v\left|t_i\right| x_i}{\mathcal{N}_v(t)}-\sum_{i \geq 1} \varphi_i\left(x_i\right) \leq v.
$$
For another situation, we may assume that $\sum_{i \geq 1} \varphi_i\left(x_i\right)=\beta v$ with $\beta>1$,
then from the properties of $\varphi(x)$ that $\varphi(0)=0$ and convexity, it follows $\sum_{i \geq 1} \varphi_i\left(x_i / \theta\right) \leq v$. In other words, $$\sum_{i \geq 1}\left|t_i\right| x_i \leq \beta \mathcal{N}_v(t),$$
so that 
$$
\sum_{i \geq 1} \frac{v\left|t_i\right| x_i}{\mathcal{N}_v(t)}-\sum_{i \geq 1} \varphi_i\left(x_i\right) \leq 0.
$$
Therefore, for any $v\geq0$, based on the definition of $\varphi^*(\cdot)$,
$$
\sum_{i \geq 1} \varphi^*_i\left(\frac{v\left|t_i\right|}{\mathcal{N}_v(t)}\right) \leq v,
$$
by choosing $\lambda=2 v /z$, we can infer that
$$
-\lambda z+\sum_{i \geq 1} \varphi^*_i\left(K \lambda\left|t_i\right|\right) \leq-2 v+v=-v,
$$
recall that $z=2 K \mathcal{N}_v(t)$, the proof is completed.

For the i.i.d. case, we observe that, from Lemma \ref{book}
$$
\begin{aligned}
\mathrm{P}\left\{Y_{t} \geq z\right\} &\leq \exp \left(\varphi^*\left( \sum_{i \geq 1} |\lambda t_i\tau_{\varphi^*}\left( X_i\right)|\right)-\lambda z\right)\\
& \leq \exp \left(\varphi^*\left( \lambda K_1 \sum_{i \geq 1} |t_i|\right)-\lambda z\right)\\
&=\exp \left(\varphi^*\left( \lambda K_1 \sum_{i \geq 1} |t_i| \right) -\lambda K_1 \sum_{i \geq 1} |t_i| \frac{z}{K_1 \sum_{i \geq 1} |t_i| }\right).
\end{aligned}
$$
Optimize over $\lambda>0$, and notice that $\varphi^*(\varphi^*(\cdot))=\varphi(\cdot)$, we obtain 
$$
\mathrm{P}\left(Y_{t} \geq z\right) \leq \exp \left(- \varphi\left(\left(\frac{z}{K_1 \left|\sum_{i \geq 1}t_i\right|}\right)\right)\right).
$$
Note that for a sufficiently small $\lambda$ such that $\lambda K_1 \leq 1$, then $\mathrm{E} \exp \lambda X_i \leq 1$. Thus, we have
$$
\begin{aligned}
\mathrm{E} \exp \lambda X_i &\leq 1+\lambda^2 \mathrm{E} X_i^2 \exp \lambda X_i \leq 1+\lambda^2\left(\mathrm{E} X_i^4\right)^{1 / 2}\left(\mathrm{E} \exp \lambda X_i\right)^{1 / 2}\\
&\leq 1+c \lambda^2 K_2^2\leq \exp c \lambda^2 K_2^2,
\end{aligned}
$$
where the first inequality is from the fact that $|\exp x-1-x| \leq |x|^2 \exp |x|$, and the second one is based on the Cauchy–Schwarz inequality. Consequently, for $0 \leq \lambda \leq \frac{1}{K_1}$,
$$
\mathrm{P}\left(Y_{t} \geq z\right) \leq \exp \left(- \frac{z^2}{K_2^2 ||t||_2^2}\right).
$$
Combine all the cases, and then we complete the proof.
\end{proof}

\begin{remark}
The quantity $\mathcal{N}_v(t)$ is introduced by Talagrand \cite{MT}, which is used to estimate the upper bound of the tail probability of the canonical process $Y_{t}=\sum_{i \geq 1} t_i X_i$. In this procedure, the log-concave property of random variable $X_i$ plays a crucial role. It is worth noting that the bound stated in Theorem \ref{2.3}, in i.i.d. case,  aligns with Bernstein's inequality for sums of independent sub-Gaussian random variables, as provided by Vershynin \cite{RV}.

It can be seen that the argument in this result for the i.i.d. case bears a resemblance to the results in the contributions of Antonini \cite{AHV},\cite{AKV}, both of which involve the use of convex conjugate transform. Moreover, it should be noted that the proof solely relies on the properties of the $N$-function, thus making little distinction between $\varphi(\cdot)$ and $\varphi^*(\cdot)$.
\end{remark}

With the tail inequality of $\mathcal{N}_v(t)$, the following result makes it possible to estimate the expected value of the supremum of the process $Y_t$. Refer to Talagrand \cite{MT1994}, \cite{MT} for more details.
\begin{corollary}\label{bv}
Let $B(v)=\left\{t ; \mathcal{N}_v(t) \leq v\right\}$, then for any $v \geq 1$ and $t \in B(v)$, we then have $||Y_t||_v \leq L u$.
\end{corollary}
\begin{proof}
By the Theorem \ref{2.3}, for any $u \geq 1$ and $w \geq v$, we have  
$$\begin{aligned}\mathrm{P}\left(Y_{t} \geq 2 u K \mathcal{N}_v(t)\right)&\leq \mathrm{P}\left(Y_{t} \geq 2 u K v\right)\\ &=\mathrm{P}\left(\frac{Y_{t}}{2u K} \geq  w\right)\leq \exp (-w),\end{aligned}$$

Denote $\frac{Y_{t}}{2u K}=\frac{Y_{t}}{2u K}1_{\{\frac{Y_{t}}{2u K} \leq v\}}+\frac{Y_{t}}{2u K}1_{\{\frac{Y_{t}}{2u K} > v\}}$. Obviously, the first part satisfies $||\frac{Y_{t}}{2u K}1_{\{\frac{Y_{t}}{2u K} \leq v\}}||_v \leq v$. For the second part, the $L_u$ norm can be derived using the tail bound $\mathrm{P}\left(\frac{Y_{t}}{2u K}1_{\{\frac{Y_{t}}{2u K} > v\}} \geq w\right) \leq \exp(-w)$ with $||\frac{Y_{t}}{2u K}1_{\{\frac{Y_{t}}{2u K} > v\}}||_u \leq L u$, so the proof is completed.
\end{proof}

Next, we consider a most special case for the $Y_{t}=\sum_{i \geq 1} t_i X_i$ with $t_i= 1$ for $i\leq N$ and $t_i=0$ for $i\geq N$, where we derive a tighter upper bound for the concentration inequality of $\sum_{i=1}^N X_i$ under $\varphi$-sub-Gaussian condition. 

\begin{theorem}\label{rad}
Let $X=(X_1, \ldots, X_n)$ be i.i.d.$\varphi$-sub-Gaussian random variables, and let $U \sim \operatorname{Unif}(0,1)$ be an independent random variable. Then, for all $\alpha \in(0,1)$, 
$$
\mathrm{P}\left( \sum_{i=1}^N (X_i-\mathrm{E}[X]) \geqslant \frac{ C \tau_{\varphi}(X) (2 \log(\frac{1}{\alpha})+\log U)}{\varphi^{(-1)}(\log(\frac{1}{\alpha}))}\right) \leqslant \alpha,
$$
where C is a constant less than $4$.
\end{theorem}

 Obviously, when $U=1$, we can retrieve the classical concentration inequality for $\varphi$-sub-Gaussian random variables, as $\log U < 0$, yielding a tighter bound. 

\begin{proof}[Proof of Theorem \ref{rad}]
According to the result in Antonini \cite{AHV}, 
we obtain that
$$
\mathrm{E} \exp \left(\lambda \sum_{i=1}^n\left(X_i-\mathrm{E}[X]\right)-\varphi( C \lambda \tau_{\varphi}(X) )\right) \leq 1
$$
for nonnegative random variables. For all $\alpha \in(0,1)$, set $\lambda = \frac{\varphi^{(-1)}(\log(\frac{1}{\alpha}))}{C\tau_{\varphi}(X)}$, then by the uniformly randomized Markov's inequality, we derive that
$$
\begin{aligned}
&\mathrm{P}\left( \sum_{i=1}^N (X_i-\mathrm{E}[X]) \geqslant \frac{ C \tau_{\varphi}(X) (2\log(\frac{1}{\alpha})+\log U)}{\varphi^{(-1)}(\log(\frac{1}{\alpha}))}\right)\\
&=\mathrm{P}\left( \sum_{i=1}^N (X_i-\mathrm{E}[X]) \geqslant  \frac{2\varphi( C \lambda \tau_{\varphi}(X) )+\log U}{\lambda}\right)\\
&=\mathrm{P}\left(\exp(\lambda \sum_{i=1}^N (X_i-\mathrm{E}[X])) \geqslant \exp( 2\varphi( C \lambda \tau_{\varphi}(X) )+\log U)\right)\\
&=\mathrm{E}[\mathrm{P}(U \leqslant \exp(\lambda \sum_{i=1}^N (X_i-\mathrm{E}[X])- 2\varphi( C \lambda \tau_{\varphi}(X) ) \mid X)]\\
&\leqslant \mathrm{E}[\exp(\lambda \sum_{i=1}^N (X_i-\mathrm{E}[X])-\varphi( C \lambda \tau_{\varphi}(X) )-\log \frac{1}{\alpha}]\\
&\leqslant \alpha.
\end{aligned}
$$
\end{proof}



 In the subsequent result, we consider concentration inequalities for the function $f: \Omega^n \rightarrow \mathbb{R}$. If $f$ is a sum of i.i.d. centered $\varphi$-sub-Gaussian variables, it can be viewed as an extension of Theorem \ref{2.3}. This result encompasses Maurer's results when $\varphi(x)=x$. To present that, we need to introduce the following Lemmas.

\begin{lemma}\label{fe}
Assume that $||X||_{\tau_{\tilde{\varphi}}}< \frac{1}{e}$ with $||X||_{\tau_{\tilde{\varphi}}}:=\sup _{p \geq 1} \frac{||X||_p}{\varphi^{(-1)}(p)}$, the funciton $\varphi(\cdot)$ satisfies $\varphi^{(-1)}(1)= 1$, and $X$ is a centered $\varphi$-sub-Gaussian variable, then 
$$
\int_0^1\left(\int_t^1 \mathrm{E}_{s X}\left[\left(X-\mathrm{E}_{s X}[X]\right)^2\right] d s\right) d t \leq \frac{e^2||X||_{\tau_{\tilde{\varphi}}}^2}{\left(1-e||X||_{\tau_{\tilde{\varphi}}}\right)^2},
$$
where $\mathrm{E}_X[Y]=\mathrm{E}\left[Y e^X\right] / \mathrm{E}\left[e^X\right]$.
\end{lemma}

\begin{proof}
Let $s \in[0,1]$, by the property of variance, we obtain that
$$
\mathrm{E}_{s X}\left[\left(X-\mathrm{E}_{s X}[X]\right)^2\right] \leq \mathrm{E}_{s X}\left[X^2\right]  .
$$

The right part of the inequality is equivalent to $\frac{\mathrm{E}\left[X^2 e^{s X}\right]}{\mathrm{E}\left[e^{s X}\right]}$, since $\mathrm{E}\left[X_k\right]=0$, by the Jensen’s inequality, $$\frac{\mathrm{E}\left[X^2 e^{s X}\right]}{\mathrm{E}\left[e^{s X}\right]}  \leq \mathrm{E}\left[X^2 e^{s X}\right]\leq \mathrm{E}\left[\sum_{m=0}^{\infty} \frac{s^m}{m !} X^{m+2}\right]=\sum_{m=0}^{\infty} \frac{s^m}{m !} \mathrm{E}\left[X^{m+2}\right]$$.

By the definition of $||X||_{\tau_{\tilde{\varphi}}}$, and $\varphi^{(-1)}(1) \leq 1 $,
$$
\begin{aligned}
\sum_{m=0}^{\infty} \frac{s^m}{m !} \mathrm{E}\left[X^{m+2}\right]& \leq \sum_{m=0}^{\infty} \frac{s^m}{m !}||X||_{\tau_{\tilde{\varphi}}}^{m+2}(\varphi^{(-1)}(m+2))^{m+2}\\
&\leq \sum_{m=0}^{\infty} \frac{s^m}{m !}||Y||_{\tau_{\tilde{\varphi}}}^{m+2}(m+2)^{m+2}\\
&\leq e^2||X||_{\tau_{\tilde{\varphi}}}^2 \sum_{m=0}^{\infty}(m+2)(m+1)\left(s e||X||_{\tau_{\tilde{\varphi}}}\right)^m.
\end{aligned}
$$
Considering the fact that  $se ||X||_{\tau_{\tilde{\varphi}}}<1$, the sum is absolute convergence. Therefore, we can change the expectation and summation, then we have 
$$
\begin{aligned}
\int_0^1\left(\int_t^1 \mathrm{E}_{s X}\left[\left(X-\mathrm{E}_{s X}[X]\right)^2\right] d s\right) d t &\leq e^2||X||_{\tau_{\tilde{\varphi}}}^2 \sum_{m=0}^{\infty}(m+1)\left(e||X||_{\tau_{\tilde{\varphi}}}\right)^m\\
&=\frac{e^2||X||_{\tau_{\tilde{\varphi}}}^2}{\left(1-e||X||_{\tau_{\tilde{\varphi}}}\right)^2},
\end{aligned}
$$
with the fact that $\int_0^1 \int_t^1 s^m d s d t=\frac{1}{m+2}$.
\end{proof}
The following Lemma is from the Maurer \cite{M2012}.
\begin{lemma}\label{m20}
Let $C_1$ and $a$ denote two positive real numbers, $t>0$. Then
$$
\inf _{\beta \in[0,1 / a)}\left(-\beta t+\frac{C_1 \beta^2}{1-a \beta}\right) \leq \frac{-t^2}{2(2 C_1+a t)} .
$$
\end{lemma}

\begin{theorem}\label{med}
Let $f: \Omega^n \rightarrow \mathbb{R}$, and assume that $X=\left(X_1, \ldots, X_n\right)$ is a vector of independent $\varphi$-sub-Gaussian random variables as in Lemma \ref{fe}. Then for any $t>0$, we have 
$$
\mathrm{P}\{f(X)-\mathrm{E} f(X^{\prime})>t\} \leq \exp \left(\frac{-t^2}{4 e^2||\sum_k|| f_k(X)||_{\tau_{\tilde{\varphi}}}^2||_{\infty}+2 e \max _k|||| f_k(X)||_{\tau_{\tilde{\varphi}}}||_{\infty} t}\right) .
$$
\end{theorem}

\begin{proof}[Proof of Theorem \ref{med}]
Our goal is to prove the upper bound of $\ln \mathrm{E}\left[e^{\beta(f-\mathrm{E} f)}\right]$ in the form of Lemma \ref{m20}, and then by the Markov inequality, we can conclude the proof.

Let $0<\zeta \leq \beta<(e \max _k|||| f_k(X)||_{\tau_{\tilde{\varphi}}}||_{\infty})^{-1}$, for any $x \in \Omega^n$ and $k \in\{1, \ldots, n\}$, we have $$||\zeta f_k(X)(x)||_{\tau_{\tilde{\varphi}}}<||f_k(X)(x)||_{\tau_{\tilde{\varphi}}} /(e\max _k|||| f_k(X)||_{\tau_{\tilde{\varphi}}}||_{\infty}) \leq 1 / e.$$

 Then we control the fluctuation representation of $\hat{h}_k(x):= \zeta f_k(X)(x)$. By the Lemma \ref{fe}, for all $x$
$$
\int_0^1\left(\int_t^1 \mathrm{E}_{s \hat{h}_k(x)}\left[\left(\hat{h}_k(x)-\mathrm{E}_{s \hat{h}_k(x)}[\hat{h}_k(x)]\right)^2\right] d s\right) d t
 \leq \frac{e^2||\hat{h}_k(x)||_{\tau_{\tilde{\varphi}}}^2}{\left(1-e||\hat{h}_k(x)||_{\tau_{\tilde{\varphi}}}\right)^2}, 
$$
then by the subadditivity introduced in Boucheron et al. \cite{BL}, 
$$
\mathrm{E}_{f(X)}[f(X)]-\ln \mathrm{E}\left[e^{f(X)}\right] \leq \mathrm{E}_{f(X)}\left[\sum_{i=1}^n \left(\mathrm{E}_{f_k(X)}[f_k(X)]-\ln \mathrm{E}\left[e^{f_k(X)}\right]\right)\right],
$$
we have
$$
\begin{aligned}
\mathrm{E}_{\zeta f(X)}[f(X)]-\ln \mathrm{E}\left[e^{\zeta f(X)}\right] &\leq \mathrm{E}_{\zeta f(X)}\left[\sum_{k} \left(\mathrm{E}_{\zeta f_k(X)}[\zeta f_k(X)]-\ln \mathrm{E}\left[e^{\zeta f_k(X)}\right]\right)\right]\\
& \leq \frac{\zeta^2 e^2 \mathrm{E}_{\zeta f(X)}\left[\sum_k||f_k\left(X\right)||_{\tau_{\tilde{\varphi}}}^2(X)\right]}{(1-\zeta e \max _k|||| f_k(X)||_{\tau_{\tilde{\varphi}}}||_{\infty})^2}\\ 
 &\leq \frac{\zeta^2 e^2||\sum_k|| f_k(X)||_{\tau_{\tilde{\varphi}}}^2||_{\infty}}{(1-\zeta e \max _k|||| f_k(X)||_{\tau_{\tilde{\varphi}}}||_{\infty})^2}.
\end{aligned}
$$
Finally, it is obtained that
$$
\ln \mathrm{E}\left[e^{\beta(f-\mathrm{E} f)}\right]=\beta \int_0^\beta \frac{S(\zeta f(X)) d \zeta}{\zeta^2} \leq \frac{\beta^2 e||\sum_k|| f_k(X)||_{\tau_{\tilde{\varphi}}}^2||_{\infty}}{1-\beta e \max _k|||| f_k(X)||_{\tau_{\tilde{\varphi}}}||_{\infty}}.
$$

\end{proof}

\section{Applications}

To present the application of our theorems, the vector-valued concentration inequality is crucial, we first consider the following lemma to derive it.
\begin{lemma}\label{2.1} Let $X, X^{*}$ be vector of independent random variables with values in space $\Omega$, $ \phi: \Omega \times \Omega \rightarrow \mathbb{R}$ measurable, then \\
(i) $||\mathrm{E}\left[\phi\left(X, X^{*}\right) \mid X\right]||_{\tau_{\varphi}} \leq||\phi\left(X, X^{*}\right)||_{\tau_{\varphi}}$.\\
(ii) If $\Omega=\mathbb{R}$ then $||X-\mathrm{E}[X]||_{\tau_{\varphi}} \leq 2||X||_{\tau_{\varphi}}$.
\end{lemma}
\begin{proof}
By the definition of $||\cdot||_{\tau_{\varphi}}$, we have 
$$
\begin{aligned}
||\mathrm{E}\left[\phi\left(X, X^{*}\right) \mid X\right]||_{\tau_{\varphi}}&= \sup _{\lambda \neq 0} \frac{\varphi^{(-1)}(\log \mathrm{E} \exp \{\lambda \mathrm{E}\left[\phi\left(X, X^{*}\right) \mid X\right]\})}{|\lambda|}\\
&\leq \sup _{\lambda \neq 0} \frac{\varphi^{(-1)}(\log \mathrm{E} \exp \{\lambda \log \mathrm{E}\left[\exp{\phi\left(X, X^{*}\right)} \mid X\right]\})}{|\lambda|}\\
&=\sup _{\lambda \neq 0} \frac{\varphi^{(-1)}(\log \mathrm{E} \exp \{\lambda\phi\left(X, X^{*}\right) \})}{|\lambda|}\\
&= ||\phi\left(X, X^{*}\right)||_{\tau_{\varphi}}.
\end{aligned}
$$

Then we assume that $\Omega=\mathbb{R}$ and $\phi(s, t)=s-t$, By using (i), we have that
$$
||X-\mathrm{E}\left[X^{*}\right]||_{\tau_{\varphi}}=||\mathrm{E}\left[X-X^{*} \mid X\right]||_{\tau_{\varphi}} \leq||X-X^{*}||_{\tau_{\varphi}} \leq 2||X||_{\tau_{\varphi}}.
$$
\end{proof}

\begin{remark}
Obviously, the norm defined in Lemma \ref{fe} also satisfies the inequality in Lemma \ref{2.1}. The proof only needs the properties of the norm and Jensen's inequality. We will not repeat it here. Interested readers can refer to \cite{MP21}.
\end{remark}

\begin{corollary}\label{7.2}
If $\mathcal{H}$ is a Hilbert space with norm $||\cdot||$, the $X_i$ are i.i.d. with $\varphi^{(-1)}(1)=1$ and $ n\geq\ln(\frac{1}{\delta})\geq 1$, then with probability at least $1-\delta$
$$
||\frac{1}{n} \sum_i X_i-\mathrm{E}\left[X_1^{*}\right]|| \leq 6 e |||| X_1||||_{\tau_{\tilde{\varphi}}} \sqrt{\frac{\ln(\frac{1}{\delta})}{n}} .
$$
where the $|||| X_1||||_{\tau_{\tilde{\varphi}}}$ is as above.
\end{corollary}
\begin{proof}
With the assumption above, by the Jensen's inequality, we have that
\begin{equation}
\begin{aligned}
\mathrm{E}\left[||\sum X_i-\mathrm{E}\left[X_i^{*}\right]||\right]&\leq \sqrt{n \mathrm{E}\left[||X_1-\mathrm{E}\left[X_i^{*}\right]||^2\right]}\\
&=\sqrt{n}||||X_1||||_2\\
&\leq \sqrt{n\ln(\mathrm{E}\left[\exp(2||\sum X_i-\mathrm{E}\left[X_i^{*}\right]||)\right])}\\
& \leq 2 \sqrt{n}||||X_1||||_{\tau_{\tilde{\varphi}}}.
\end{aligned}
\end{equation}

Let $f(x)=||\sum_i\left(x_i-\mathrm{E}\left[X_1^{*}\right]\right)||$, then
\begin{equation}
\begin{aligned}
\left|f_k(X)(x)\right|&=|||\sum_{i \neq k} x_i+X_k-n \mathrm{E}\left[X_1^{*}\right]||-\mathrm{E}[||\sum_{i \neq k} x_i+X_k^{\prime}-n \mathrm{E}\left[X_1^{*}\right]||]| \\
&\leq \mathrm{E}\left[||X_k-X_k^{*}|| \mid X\right].
\end{aligned}
\end{equation}

According to the properties of the norm, by the Lemma \ref{2.1} and the Theorem \ref{med}, with probability at least $1-\delta$, we have
\begin{equation}
\begin{aligned}
||\sum_i X_i-\mathrm{E}\left[X_1^{*}\right]|| & \leq \mathrm{E}\left[||\sum X_i-\mathrm{E}\left[X_i^{*}\right]||\right]+2 e|||| X_1||||_{\tau_{\tilde{\varphi}}} \sqrt{n \ln (1 / \delta)}+2 e|||| X_1||||_{\tau_{\tilde{\varphi}}} \ln (1 / \delta) \\
& \leq 2 \sqrt{n}||||X_1||||_{\tau_{\tilde{\varphi}}}+2 e|||| X_1||||_{\tau_{\tilde{\varphi}}} \sqrt{n \ln (1 / \delta)}+2 e|||| X_1||||_{\tau_{\tilde{\varphi}}} \ln (1 / \delta) \\
& \leq 2 \sqrt{n}|||| X_1||||_{\tau_{\tilde{\varphi}}}(1+2e\sqrt{ \ln (1 / \delta)})) \\
& \leq 6 e |||| X_1||||_{\tau_{\tilde{\varphi}}}  \sqrt{\frac{\ln(\frac{1}{\delta})}{n}},
\end{aligned}
\end{equation}
where the last inequality is from $\ln(\frac{1}{\delta})\geq \frac{1}{e}$, finally, divide both sides by $n$, we complete the proof.
\end{proof}


\subsection{Uniform bound for PCA}
Based on the results presented, we now prepare for deriving an upper bound for principal subspace selection, commonly referred to as PCA, particularly when dealing with the data following a $\varphi$-sub-Gaussian distribution. PCA involves finding a projection onto a $d$-dimensional subspace within a given $n$-dimensional data space, allowing the subspace to capture the properties of the data. Denote $\mathcal{H}$ as a Hilbert space, where $X_i$ are independent and identically distributed random variables taking values in $\mathcal{H}$. We define $\mathcal{P}_d$ as the set of $d$-dimensional orthogonal projection operators within $\mathcal{H}$. For any $x \in \mathcal{H}$ and $\mathrm{P} \in \mathcal{P}_d$, the reconstruction error is quantified using a loss function $\ell(\mathrm{P}, x):=||\mathrm{P} x-x||_{\mathcal{H}}^2$. In the following result, we establish an upper bound on the disparity between the expected and empirical reconstruction errors, ensuring uniformity across all projections within $\mathcal{P}_d$.

\begin{theorem} Let $X=\left(X_1, \ldots, X_n\right)$ be i.i.d. with $\varphi^{(-1)}(1)= 1$ and $n \geq \ln (1 / \delta) \geq 1$, then with probability at least $1-\delta$, we have
$$
\begin{aligned}
\sup _{\mathrm{P} \in \mathcal{P}_d} \frac{1}{n} \sum_i \ell\left(\mathrm{P}, X_i\right) - \mathrm{E}\left[\ell\left(\mathrm{P}, X_1\right)\right] \leq 12 \sqrt{d}  e K_3 \sqrt{\frac{\ln(\frac{1}{\delta})}{n}},
\end{aligned}
$$

where $K_3=\max_{ i} \tau_{\tilde{\varphi}}(|| X_1||^2)$. 
\end{theorem}
\begin{proof}
We consider the space of Hilbert-Schmidt operators, allowing us to express $\ell(\mathrm{P}, x)$ as $||R_x||_{H S}-\left\langle \mathrm{P}, R_x\right\rangle_{H S}$, then we have
\begin{equation}\label{3.1}
\begin{array}{rl}
\sup _{\mathrm{P} \in \mathcal{P}_d} \frac{1}{n}  \sum_i \ell\left(\mathrm{P}, X_i\right) - \mathbb{E} \left[\ell\left(\mathrm{P}, X_1\right)\right]
 &= \sup _{\mathrm{P} \in \mathcal{P}_d}\left\langle \mathrm{P}, \frac{1}{n} \sum_i\left(R_{X_i}-\mathbb{E}\left[R_{X_i}\right]\right)\right\rangle_{H S}\\
&+\left(\mathbb{E}||R_{X_i}||_{H S}-\frac{1}{n} \sum_i||R_{X_i}||_{H S}\right) .
\end{array}
\end{equation}

Since $||\mathrm{P}||_{H S}=\sqrt{d}$ for $\mathrm{P} \in \mathcal{P}_d$, we can apply the Cauchy-Schwarz inequality and the Corollary \ref{7.2} to bound the first term of (\ref{3.1}) by 
$$
\sqrt{d}||\frac{1}{n} \sum_i\left(R_{X_i}-\mathrm{E}\left[R_{X_1}\right]\right)||_{H S} \leq \sqrt{d} 6 e ||||R_{X_i}||_{H S}||_{\tau_{\tilde{\varphi}}} \sqrt{\frac{\ln(\frac{1}{\delta})}{n}},
$$
with probability at least $1-\delta$. Using a similar argument as in Corollary \ref{7.2}, the second term can by  bounded by $||R_{X_i}||_{H S}$ in the Hilbert space, then by applying Theorem \ref{2.3} to a sum of independent $\varphi$-sub-Gaussian, note that $ ||||R_{X_1}||_{H S}||_{\tau_{\tilde{\varphi}}}=|||| X_1||^2||_{\tau_{\tilde{\varphi}}}$, since $\ln (1 / \delta) \geq 1$, we can infer that $\varphi^{(-1)}(\ln (1 / \delta)) \leq \ln (1 / \delta)$. Thus, we obtain,
$$
\begin{aligned}
\left(\mathbb{E}||Q_{X_i}||_{H S}-\frac{1}{n} \sum_i||Q_{X_i}||_{H S}\right)&\leq \max\left(K_3 \sqrt{\frac{\log(\frac{1}{\delta})}{n}},\max _i|||| X_1||^2||_{\psi_1} \sqrt{\frac{\log(\frac{1}{\delta})}{n}}\right) \\
&\leq K_3 \sqrt{\frac{\log(\frac{1}{\delta})}{n}},
\end{aligned}$$
the result follows from a union bound with probability at least $1-\delta$.
\end{proof}

\subsection{Application in Rademacher complexities}
Let $\mathcal{G}$ be a class of functions $g: \mathcal{G} \rightarrow \mathbb{R}$. We seek a bound on the supremum deviation of the following random variable, with a high probability
$$
f(X)=\sup _{g \in \mathcal{G}} \frac{1}{n} \sum_{i=1}^n\left(g\left(X_i\right)-\mathrm{E}\left[g\left(X_i^{*}\right)\right]\right) .
$$

The classical method of Rademacher complexities requires expressing $f(X)$ in the following form
$$
f(X)=\left(f(X)-\mathrm{E}\left[f\left(X^{*}\right)\right]\right)+\mathrm{E}\left[f\left(X^{*}\right)\right].
$$
The first term mentioned above can be bounded by the concentration inequality. In contrast to the classical approach using the bounded difference inequality, we only require Lipschitz properties of $g(X)$ for $\varphi$-sub-Gaussian random variables. The second term, $\mathrm{E}[f(X)]$, arises from the symmetrization argument. For the independent Rademacher variables $\epsilon_i$, uniformly distributed on $\{-1,1\}$, we have
$$
\mathrm{E}[f(X)] \leq \mathrm{E}\left[\frac{2}{n} \mathrm{E}\left[\sup _{g \in \mathcal{G}} \sum_i \epsilon_i g\left(X_i\right) \mid X\right]\right].
$$
Denote the quantity above as $\mathrm{E}[\mathcal{L}(\mathcal{G}, X)]$, it can be observed that $\mathrm{E}[\mathcal{L}(\mathcal{G}, X)]$ varies as the classes $\mathcal{G}$ change.

\begin{theorem}\label{th9} Let $X=\left(X_1, \ldots, X_n\right)$ be random vector as above with values in a Banach space $(\Omega,||\cdot||)$ and let $\mathcal{G}$ be a class of Lipschitz functions $g: \Omega \rightarrow \mathbb{R}$, i.e. $g(x)-g(y) \leq L||x-y||$ for all $g \in \mathcal{G}$ and $x, y \in \Omega$. If $n \geq \ln (1 / \delta)\geq1$ then with probability at least $1-\delta$
$$
\sup _{g \in \mathcal{G}} \frac{1}{n} \sum_i g\left(X_i\right)-\mathrm{E}(g(X)) \leq \mathrm{E}[\mathcal{L}(\mathcal{G}, X)]+12 e L |||| X_1||||_{\tau_{\tilde{\varphi}}}  \sqrt{\frac{\ln(\frac{1}{\delta})}{n}}.
$$
\end{theorem}

\begin{proof}
The vector space
$$
\mathcal{S}=\left\{h: \mathcal{G} \rightarrow \mathbb{R}: \sup _{g \in \mathcal{G}}|h(g)|<\infty\right\}
$$
becomes a normed space with norm $||h||_{\mathcal{S}}=\sup _{g \in \mathcal{G}}|h(g)|$. For each $X_i$, define $\hat{X}_i \in \mathcal{S}$ by $\hat{X}_i(g)=(1 / n)\left(g\left(X_i\right)-\mathrm{E}\left[g\left(X_i^{\prime}\right)\right]\right)$. Then the $\hat{X}_i$ are zero mean random variables in $\mathcal{S}$ and $f(X)=||\sum_i \hat{X}_i||_{\mathcal{S}}$. Also with Lemma \ref{2.1} and the i.i.d.-assumption,

$$
\begin{aligned}
|||| \hat{X}_i||_{\mathcal{S}}||_{\tau_{\tilde{\varphi}}} & =\frac{1}{n}||\sup _h\left(\mathrm{E}\left[g\left(X_i\right)-g\left(X_i^{*}\right)\right] \mid X\right)||_{\tau_{\tilde{\varphi}}} \\
& \leq \frac{L}{n}||\mathrm{E}\left[||X_i-X_i^{\prime}||\right] \mid X||_{\tau_{\tilde{\varphi}}} \leq \frac{2 L}{n}|||| X_1||||_{\tau_{\tilde{\varphi}}}.
\end{aligned}
$$

From Corollary \ref{7.2}, we can deduce that with probability at least $1-\delta$, the following holds,
$$
f(X)-\mathrm{E}\left[f\left(X^{*}\right)\right] \leq 6 e L |||| X_1||||_{\tau_{\tilde{\varphi}}} \sqrt{\frac{\ln(\frac{1}{\delta})}{n}}.
$$
\end{proof}

Specifically, we consider the case of linear regression with unbounded data. Let $\Omega=(\mathcal{G}, \mathbb{R})$, where $\mathcal{G}$ is defined as above with inner product $\langle.,.\rangle$ and norm $||.||_{\mathcal{G}}$. Additionally, let $X_1$ and $Y_1$ be $\varphi$-sub-Gaussian random variables in $\mathcal{G}$ and $\mathbb{R}$. The pair $\left(X_1, Y_1\right)$ represents the joint of input-vectors $X_1$ and outputs $Y_1$. $\mathcal{G}=\left\{(x, y) \mapsto g(x, y)=\ell(\langle w, x\rangle-y):||w||_{\mathcal{G}} \leq L\right\}$ is the class of functions defined on $\Omega$, where $\ell$ is a 1-Lipschitz loss function.

\begin{corollary}
 Let $\Omega$ and $\mathcal{G}$ be defined as above and $(X, Y)=\left(\left(X_1, Y_1\right), \ldots,\left(X_n, Y_n\right)\right)$ be an i.i.d. sample of random variables in $\Omega$. Then for $\delta>0$ and $n \geq \ln (1 / \delta)\geq 1$ with probability at least $1-\delta$
$$
\sup _{g \in \mathcal{G}} \frac{1}{n} \sum_i g\left(X_i, Y_i\right)-\mathrm{E}\left(g\left(X_i, Y_i\right)\right) \leq \frac{12}{\sqrt{n}}\left(L|||| X_1||||_{\tau_{\tilde{\varphi}}}+|||| Y_1||||_{\tau_{\tilde{\varphi}}}\right)(1+e \sqrt{\ln(\frac{1}{\delta})}).
$$
\end{corollary}

\begin{proof}
Considering $\Omega$ as a Banach space with the norm $||(x, y)||=L||x||_H+|y|$, then we have $||||\left(X_1, Y_1\right)||||_{\tau_{\tilde{\varphi}}} \leq$ $L|||| X_1||||_{\tau_{\tilde{\varphi}}}+||\left|Y_1\right|||_{\tau_{\tilde{\varphi}}}$ which is based on the properties of the norm, for any $g \in \mathcal{G}$
$$
\begin{aligned}
g(x, y)-g\left(x^{*}, y^{*}\right) & =\ell(\langle w, x\rangle-y)-\ell\left(\left\langle w, x^{*}\right\rangle-y^{*}\right) \\
& \leq L||x-x^{*}||_{\mathcal{G}}+\left|y-y^{*}\right| \leq||(x, y)-\left(x^{*}, y^{*}\right)||,
\end{aligned}
$$

so the class $\mathcal{G}$ is uniformly 1-Lipschitz. For an i.i.d. sample $(X, Y) \in \Omega^n$, using the Lipschitz property of $\ell$ and the Jensen's inequality, we obtain that
$$
\mathcal{L}(\mathcal{G},(X, Y)) \leq \frac{2}{n}\left(L \sqrt{\sum_i||X_i||_{\mathcal{G}}^2}+\sqrt{\sum_i\left|Y_i\right|^2}\right) .
$$

Through a similar argument as in the Corollary \ref{7.2},
$$
\mathrm{E}[\mathcal{L}(\mathcal{G},(X, Y))] \leq \frac{12}{\sqrt{n}}\left(L||||| X_1||||_{\tau_{\tilde{\varphi}}}+|||| Y_1||||_{\tau_{\tilde{\varphi}}}\right) .
$$

Then by the Throrem \ref{th9}, with probability at least $1-\delta$,
$$
\sup _{g \in \mathcal{G}} \frac{1}{n} \sum_i g\left(X_i, Y_i\right)-\mathrm{E}\left(g\left(X_i, Y_i\right)\right) \leq \frac{12}{\sqrt{n}}\left(L|||| X_1||||_{\tau_{\tilde{\varphi}}}+|||| Y_1||||_{\tau_{\tilde{\varphi}}}\right)(1+e \sqrt{\ln(\frac{1}{\delta})}) .
$$
\end{proof}

This result can be extended to other classic problems involving Rademacher complexity in the empirical process, such as classification and clustering problems. Additionally, we can apply this result by replacing Rademacher complexity with VC-dimension complexity as in the work of Depeisin \cite{dr}.

\section *{Acknowledgments}
The authors would like to express their sincere gratitude to Prof. Hanchao Wang for his suggestions and fruitful discussions.

\end{document}